\documentclass[12pt]{amsart}

\copyrightinfo{2000}{American Mathematical Society}

\newcommand{\cD}{{\mathcal D}}
\newcommand{\cE}{{\mathcal E}}
\newcommand{\cL}{{\mathcal L}}

\newcommand{\bN}{{\mathbb N}}
\newcommand{\bZ}{{\mathbb Z}}
\newcommand{\bQ}{{\mathbb Q}}
\newcommand{\bR}{{\mathbb R}}
\newcommand{\bC}{{\mathbb C}}

\numberwithin{equation}{section}

\newtheorem{Theorem}{Theorem}[section]
\newtheorem{Lemma}{Lemma}[section]

\newtheorem{Definition}{Definition}[section]
\newtheorem{Remark}{Remark}[section]

\newtheorem{Proposition}{Proposition}[section]

\author{V.~M.~Shelkovich}
\address{Department of Mathematics, St.-Petersburg State Architecture
and Civil Engineering University, \ 2 Krasnoarmeiskaya 4, 190005,
St. Petersburg, \ Russia. \ Phone: +7\,(812)\,2517549 \,
Fax: +7\,(812)\,3165872}
\email{shelkv@vs1567.spb.edu}

\author{M.~Skopina}
\address{Department of Applied Mathematics and Control Processes,
St. Petersburg State University, \ Universitetskii pr.-35,
Petrodvorets, 198504 St. Petersburg, Russia. \ Phone: +7\,(812)\,51326090 \,
Fax: +7\,(812)\,}
\email{skopina@ms1167.spb.edu}

\title[$p$-Adic Haar multiresolution analysis]
{$p$-Adic Haar multiresolution analysis}

\thanks{The first author (V.~S.) was also supported in part by
DFG Project 436 RUS 113/809 and Grant 05-01-04002-NNIOa of Russian
Foundation for Basic Research.}



\date{ }

\begin{document}

\begin{abstract}
In this paper, the notion of {\em  $p$-adic multiresolution
analysis (MRA)} is  introduced. We use a ``natural'' refinement
equation whose solution (a refinable function) is the
characteristic function of the unit disc. This equation reflects
the fact that the characteristic function of the unit disc is
the sum of $p$  characteristic  functions of
disjoint discs of radius $p^{-1}$. The case $p=2$ is studied in
detail. Our MRA is a $2$-adic analog of the real Haar MRA. But in
contrast to the real setting, the refinable function generating
our Haar MRA is periodic with  period $1$, which never holds
for real refinable functions. This fact implies that there exist
infinity many different $2$-adic orthonormal wavelet bases in
${\cL}^2(\bQ_2)$ generated by the same Haar MRA. All of these bases
are constructed. Since  $p$-adic pseudo-differential operators
are closely related to wavelet-type bases,
our bases can be intensively used for applications.
\end{abstract}

\maketitle

\section{Introduction}
\label{s1}

\subsection{$p$-Adic wavelets and pseudo-differential operators.}\label{s1.1}
According to the well-known Ostrovsky theorem, {\it any nontrivial
valuation on the field ${\bQ}$ is equivalent either to the real
valuation $|\cdot|$ or to one of the $p$-adic valuations $|\cdot|_p$}.
We recall that the field $\bQ_p$ of $p$-adic numbers is defined as
the completion of the field of rational numbers $\bQ$ with respect
to the non-Archimedean $p$-adic norm $|\cdot|_p$. This norm is
defined as follows: if an arbitrary rational number $x\ne 0$ is
represented as $x=p^{\gamma}\frac{m}{n}$, where $\gamma=\gamma(x)\in \bZ$,
and $m$ and $n$ are not divisible by $p$, then
\begin{equation}
\label{1}
|x|_p=p^{-\gamma}, \quad x\ne 0, \qquad |0|_p=0.
\end{equation}
This norm in $\bQ_p$ satisfies the strong triangle inequality
$|x+y|_p\le \max(|x|_p,|y|_p)$.

Thus there are two equal in rights universes: the real universes and
the $p$-adic one. The latter has a specific and unusual properties.
Nevertheless, there are a lot of papers where different applications
of $p$-adic analysis to physical problems, stochastics, cognitive
sciences and psychology are studied~\cite{Ar-Dr-V}--~\cite{Bik-V},
~\cite{Kh1}--~\cite{Koch3},~\cite{Vl-V-Z}--~\cite{V2} (see also the
references therein). In view of the Ostrovsky theorem such investigations
not only have great interest in itself, but lead to applications and
better understanding of similar problems in {\it usual} mathematical
physics.

We recall that there exists a $p$-adic analysis connected with
the mapping $\bQ_p$ into $\bQ_p$ and an analysis connected with the
mapping $\bQ_p$ into the field of complex numbers $\bC$, there exist
two types of $p$-adic physics models.
For the $p$-adic analysis related to the mapping $\bQ_p \to \bC$ the
operation of partial differentiation is {\it not defined\/}, and as
a result, large number of models connected with $p$-adic differential
equations use pseudo-differential operators and the theory of
$p$-adic distributions (generalized functions) (see the above
mentioned papers and books). In particular, fractional operators
$D^{\alpha}$ are extensively used in applications (see fore-quoted
papers and especially~\cite{Vl-V-Z}).

It is well known that the theory of $p$-adic pseudo-differential operators
(in particular, fractional operators) and equations closely related to
wavelet type bases. It is typical that $p$-adic compactly
supported wavelets are eigenfunctions of $p$-adic pseudo-differential
operators~\cite{Al-Kh-Sh3}--~\cite{Al-Kh-Sh5},~\cite{Kh-Koz1},
~\cite{Kh-Koz2},~\cite{Kh-Sh1},~\cite{Koz0} --~\cite{Koz2}.
Thus the wavelet theory plays a key role in application of $p$-adic
analysis and gives a new powerful technique for solving $p$-adic
problems. This theory starts development only in resent years
and has many open problems.

In~\cite{Koz0}, S.~V.~Kozyrev constructed the orthonormal compactly
supported $p$-adic wavelet basis (\ref{62.0-1}) in ${\cL}^2(\bQ_p)$:
\begin{equation}
\label{62.0-1}
\theta_{\gamma j a}(x)=p^{-\gamma/2}\chi_p\big(p^{-1}j(p^{\gamma}x-a)\big)
\Omega\big(|p^{\gamma}x-a|_p\big), \quad x\in \bQ_p,
\end{equation}
$j\in J_{p}=\{1,2,\dots,p-1\}$, $\gamma\in \bZ$, $a\in I_p=\bQ_p/\bZ_p$.
Kozyrev's wavelets (\ref{62.0-1}) are eigenfunctions of the Vladimirov
fractional operator~\cite[IX]{Vl-V-Z}. Further development and generalization
of the theory of such type wavelets can be found in the papers by
S.~V.~Kozyrev~\cite{Koz1},~\cite{Koz2}, A.~Yu.~Khrennikov, and
S.~V.~Kozyrev \cite{Kh-Koz1}, ~\cite{Kh-Koz2}, J.~J.~Benedetto, and
R.~L.~Benedetto~\cite{Ben-Ben}, and R.~L.~Benedetto \cite{Ben1}.
In~\cite{Al-Kh-Sh3}, the multidimensional $p$-adic wavelets
generated by direct product of the Kozyrev one-dimensional
wavelets were introduced. In~\cite{Kh-Sh1}, a new type of $p$-adic
multidimensional wavelet basis was introduced:
$$
\theta_{\gamma s a}^{(m)}(x)=p^{-\gamma/2}\chi_p\big(s(p^{\gamma}x-a)\big)
\Omega\big(|p^{\gamma}x-a|_p\big), \quad x\in \bQ_p,
$$
where $s\in J_{p;m}$, $\gamma\in \bZ$, $a\in I_p$.
Here
$J_{p;m}=\{s=p^{-m}\big(s_{0}+s_{1}p+\cdots+s_{m-1}p^{m-1}\big):
s_j=0,1,\dots,p-1; \, j=0,1,\dots,m-1; s_0\ne 0\}$, \
$m\ge 1$ is a {\it fixed} positive integer.
The multidimensional wavelets from~\cite{Al-Kh-Sh3} are a particular
case of the last wavelets.
Moreover, in~\cite{Al-Kh-Sh3},~\cite{Kh-Sh1}, there were derived the
necessary and sufficient conditions for a class of multidimensional
$p$-adic pseudo-differential operators (including fractional operator)
to have such multidimensional wavelets as eigenfunctions.

It remains to point out that for pseudo-differential operators
from~\cite{Al-Kh-Sh3},~\cite{Kh-Sh1} a ``natural'' definition domain
is the Lizorkin spaces of distributions $\Phi'(\bQ_p^n)$, introduced
in~\cite{Al-Kh-Sh3}. The space $\Phi'(\bQ_p^n)$ is
{\it invariant\/} under the mentioned above pseudo-differential operators.
Moreover, the above mentioned $p$-adic wavelets
belong to the Lizorkin space $\Phi(\bQ_p^n)$ of test functions.
Recall that the {\it usual} Lizorkin spaces were studied in the
excellent papers of P.~I.~Lizorkin~\cite{Liz1},~\cite{Liz3}
(see also~\cite{Sam3},~\cite{Sam-Kil-Mar}).

It's interesting to compare appearing  first wavelets in $p$-adic
analysis with the history of the wavelet theory in  real analysis.
In 1910 Haar~\cite{25} constructed an orthogonal basis for
$\cL_2(\bR)$ consisting of the dyadic shifts and scales of one
piecewise constant function. A lot of mathematicians actively
studied Haar basis, different kinds of generalizations were
introduced, but during almost the whole century nobody could find
another wavelet function (a function whose shifts and scales form
an orthogonal basis). Only in early nineties a method for
construction of wavelet functions appeared. This method is based
on the notion of multiresolution analysis (MRA in the sequel)
introduced by Y.~Meyer and S.~Mallat~\cite{13}, \cite{15},
\cite{18}. Smooth compactly supported wavelet functions were found
in this way, which has been very important for some engineering
applications. In this paper we introduce MRA in $\cL_2(\bQ_p)$ and
present a concrete MRA for $p=2$ being an analog of Haar MRA in
$\cL_2(\bR)$. The same scheme as  in the real setting leads to a
Haar basis. It turned out that this Haar basis coincides with
Kozyrev's wavelet system. However, 2-adic Haar MRA is not an
identical copy of its real analog. In contrast to Haar MRA in
$\cL_2({\bR})$, we proved that there exist {\it infinity many
different Haar orthogonal bases} in $\cL_2({\bQ}_2)$ generated by
the same MRA.

\subsection{Contents of the paper.}\label{s1.2}
In Sec.~\ref{s2}, we recall some facts from the $p$-adic theory of
distributions~\cite{G-Gr-P}, \cite{Taib1}, \cite{Taib3}, \cite{Vl-V-Z}.
In Sec.~\ref{s3}, some facts from the theory of the $p$-adic Lizorkin
spaces~\cite{Al-Kh-Sh3} are recalled.

In Sec.~\ref{s4}, by Definition~\ref{de1} we introduce the MRA adapted
to the $p$-adic case.
In Subsec.~\ref{s4.2}, we introduce the {\it refinement equation}
(\ref{62.0-3})
$$
\phi(x)=\sum_{r=0}^{p-1}\phi\Big(\frac{1}{p}x-\frac{r}{p}\Big),
\quad x\in \bQ_p,
$$
whose solution $\phi(x)=\Omega\big(|x|_p\big)$ is the
characteristic function of the unit disc, where where $\Omega(t)$
is the characteristic function of the interval $[0,1]$.
The conjecture to use the above equation as the {\it refinement equation}
was proposed in~\cite{Kh-Sh1}.
The above {\it refinement equation} is {\it natural} and reflects
the fact that the characteristic function $\Omega\big(|x|_p\big)$
of the unit disc $B_{0}$ is represented as a sum of $p$ pieces
characteristic functions of the disjoint discs $B_{-1}(r)$,
$r=0,1,\dots,p-1$ (see (\ref{79})).

In Subsec.~\ref{s4.3}, the $2$-adic MRA is constructed.
Namely, we proved that MRA is generated by a {\it refinable function}
which is the characteristic function $\phi(x)=\Omega\big(|x|_2\big)$
of the unit disc $B_{0}=\{x: |x|_2 \le 1\}\subset {\bQ}_2$ and satisfies
the {\it refinement equation} (\ref{62.0-4})
$$
\phi(x)=\sum_{r=0}^{1}\phi\Big(\frac{1}{2}x-\frac{r}{2}\Big),
\quad x\in {\bQ}_2.
$$
By our MRA we construct $2$-adic orthonormal wavelet basis (\ref{62.0-7})
in ${\cL}^2(\bQ_2)$, which is the Kozyrev basis (\ref{62.0-1}) for the
case $p=2$.
It turned out that the Kozyrev wavelet basis is not unique orthonormal
wavelet basis.

In Sec.~\ref{s5}, {\it infinity many different $2$-adic wavelet orthonormal
bases} in ${\cL}^2(\bQ_2)$ are constructed.
Namely, using Theorem~\ref{th4}, we construct wavelet functions $\psi^{(s)}(x)$,
$s\in \bN$ whose dilatations and shifts form $2$-adic orthonormal wavelet
bases in ${\cL}^2(\bQ_2)$.

Since many $p$-adic models use pseudo-differential operators,
in particular, fractional operator, these results on $p$-adic wavelets can be
intensively used in applications. Moreover, $p$-adic wavelets can be used to
construct solutions of linear and {\it semi-linear} pseudo-differential
equations~\cite{Al-Kh-Sh5},~\cite{Koz-Os-Av-1}.

\section{$p$-Adic distributions}
\label{s2}

We recall some facts from the theory of $p$-adic distributions
(generalized functions). Here and in what follows, we shall
systematically use the notations and results from~\cite{Vl-V-Z}
and~\cite[Ch.II]{G-Gr-P}.
Let $\bN$, $\bZ$, $\bC$ be the sets of positive integers, integers,
complex numbers, respectively, and $\bN_0=\{0\}\cup\bN$.
Denote by $\bQ_p^{*}=\bQ_p\setminus\{0\}$ the multiplicative group
of the field $\bQ_p$.

The canonical form of a $p$-adic number $x\ne 0$ is
\begin{equation}
\label{2}
x = p^{\gamma}(x_0 + x_1p + x_2p^2 + \cdots),
\end{equation}
where $\gamma=\gamma(x)\in \bZ$, \ $x_j=0,1,\dots,p-1$, $x_0\ne 0$,
$j=0,1,\dots$. The series is convergent in the $p$-adic norm (\ref{1}),
and one has $|x|_p=p^{-\gamma}$.
By means of representation (\ref{2}), the {\it fractional part} $\{x\}_p$
of a number $x\in \bQ_p$ is defined as follows
\begin{equation}
\label{8.2**}
\{x\}_p=\left\{
\begin{array}{lll}
0,\quad \text{if} \quad \gamma(x)\ge 0 \quad  \text{or} \quad x=0,&&  \\
p^{\gamma}(x_0+x_1p+x_2p^2+\cdots+x_{|\gamma|-1}p^{|\gamma|-1}),
\quad \text{if} \quad \gamma(x)<0. && \\
\end{array}
\right.
\end{equation}

The function
\begin{equation}
\label{8.2-1**}
\chi_p(\xi x)=e^{2\pi i\{\xi x\}_p}
\end{equation}
for every fixed $\xi \in \bQ_p$ is an {\it additive character} of
the field $\bQ_p$.

According to~\cite[III.2.]{Vl-V-Z}, any {\it multiplicative character}
$\pi$ of the field $\bQ_p$ can be represented as
$$
\pi(x)\stackrel{def}{=}\pi_{\alpha}(x)=|x|_p^{\alpha-1}\pi_{1}(x),
\quad x \in \bQ_p^{*},
$$
where $\pi(p)=p^{1-\alpha}$ and $\pi_{1}(x)$ is a
{\it normed multiplicative character\/} such that
$\pi_1(x)=\pi_{1}(|x|_px)$, $\pi_1(p)=\pi_1(1)=1$, $|\pi_1(x)|=1$.
We denote $\pi_{0}=|x|_p^{-1}$.

The space $\bQ_p^n=\bQ_p\times\cdots\times\bQ_p$ consists of points
$x=(x_1,\dots,x_n)$, where $x_j \in \bQ_p$, $j=1,2\dots,n$, \ $n\ge 2$.
The $p$-adic norm on $\bQ_p^n$ is
\begin{equation}
\label{8}
|x|_p=\max_{1 \le j \le n}|x_j|_p, \quad x\in \bQ_p^n,
\end{equation}
where $|x_j|_p$ id defined by (\ref{1}).

Denote by $B_{\gamma}^n(a)=\{x\in \bQ_p^n: |x-a|_p \le p^{\gamma}\}$
the ball of radius $p^{\gamma}$ with the center at a point
$a=(a_1,\dots,a_n)\in \bQ_p^n$ and by
$S_{\gamma}^n(a)=\{x\in \bQ_p^n: |x-a|_p = p^{\gamma}\}
=B_{\gamma}^n(a)\setminus B_{\gamma-1}^n(a)$ its boundary (sphere),
$\gamma \in \bZ$. For $a=0$ we set $B_{\gamma}^n(0)=B_{\gamma}^n$ and
$S_{\gamma}^n(0)=S_{\gamma}^n$. For the case $n=1$ we will omit the
upper index $n$.
It is clear that
\begin{equation}
\label{9}
B_{\gamma}^n(a)=B_{\gamma}(a_1)\times\cdots\times B_{\gamma}(a_n),
\end{equation}
where $B_{\gamma}(a_j)=\{x_j: |x_j-a_j|_p \le p^{\gamma}\}\subset\bQ_p$
is a disc of radius $p^{\gamma}$ with the center at a point $a_j\in \bQ_p$,
$j=1,2\dots,n$.

Any two balls in $\bQ_p^n$ either are disjoint or one
contains the other. Every point of the ball is its center.

According to~\cite[I.3,Examples 1,2.]{Vl-V-Z}, the disc $B_{\gamma}$
is represented by the sum of $p^{\gamma-\gamma'}$ {\it disjoint}
discs $B_{\gamma'}(a)$, $\gamma'<\gamma$:
\begin{equation}
\label{79.0}
B_{\gamma}=B_{\gamma'}\cup\cup_{a}B_{\gamma'}(a),
\end{equation}
where $a=0$ and
$a=a_{-r}p^{-r}+a_{-r+1}p^{-r+1}+\cdots+a_{-\gamma'-1}p^{-\gamma'-1}$
are the centers of the discs $B_{\gamma'}(a)$, \
$r=\gamma,\gamma-1,\gamma-2,\dots,\gamma'+1$, \, $0\le a_j\le p-1$,
\, $a_{-r}\ne 0$.
In particular, the disc $B_{0}$ is represented by the sum of $p$
{\it disjoint} discs
\begin{equation}
\label{79}
B_{0}=B_{-1}\cup\cup_{r=1}^{p-1}B_{-1}(r),
\end{equation}
where $B_{-1}(r)=\{x\in S_{0}: x_0=r\}=r+p\bZ_p$, $r=1,\dots,p-1$;
$B_{-1}=\{|x|_p\le p^{-1}\}=p\bZ_p$; and
$S_{0}=\{|x|_p=1\}=\cup_{r=1}^{p-1}B_{-1}(r)$. Here all the discs
are disjoint. We call coverings (\ref{79.0}) and (\ref{79}) the
{\it canonical covering} of the discs $B_{0}$ and $B_{\gamma}$,
respectively.

On $\bQ_p$ there exists the Haar measure, i.e., a positive measure $dx$
invariant under shifts, $d(x+a)=dx$, and normalized by the equality
$\int_{|\xi|_p\le 1}\,dx=1$.
The invariant measure $dx$ on the field $\bQ_p$ is extended to an
invariant measure $d^n x=dx_1\cdots dx_n$ on $\bQ_p^n$ in the standard way.

If $f$ is an integrable function on $\bQ_p$, then~\cite[Ch.II,\S 2.2]{G-Gr-P},
~\cite[IV]{Vl-V-Z}:
\begin{equation}
\label{9.1-1}
\begin{array}{rcl}
\displaystyle
\int_{B_{\gamma}}\,dx&=&p^\gamma, \\
\displaystyle
\int_{B_{N}}f(x)\,dx&=&\sum_{\gamma=-\infty}^{N}\int_{S_{\gamma}}f(x)\,dx, \\
\displaystyle
\int_{S_{\gamma}}f(x)\,dx&=&
\int_{B_{\gamma}}f(x)\,dx-\int_{B_{\gamma-1}}f(x)\,dx. \\
\end{array}
\end{equation}

A complex-valued function $f$ defined on $\bQ_p^n$ is called
{\it locally-constant} if for any $x\in \bQ_p^n$ there exists
an integer $l(x)\in \bZ$ such that
$$
f(x+y)=f(x), \quad y\in B_{l(x)}^n.
$$

Let ${\cE}(\bQ_p^n)$ and ${\cD}(\bQ_p^n)$ be the linear spaces of
locally-constant $\bC$-valued functions on $\bQ_p^n$ and locally-constant
$\bC$-valued functions with compact supports (so-called test functions),
respectively~\cite[VI.1.,2.]{Vl-V-Z}. If $\varphi \in {\cD}(\bQ_p^n)$,
according to Lemma~1 from~\cite[VI.1.]{Vl-V-Z}, there exists $l\in \bZ$,
such that
$$
\varphi(x+y)=\varphi(x), \quad y\in B_l^n, \quad x\in \bQ_p^n.
$$
The largest of such numbers $l=l(\varphi)$ is called the
{\it parameter of constancy} of the function $\varphi$.
Let us denote by ${\cD}^l_N(\bQ_p^n)$ the finite-dimensional space of
test functions from ${\cD}(\bQ_p^n)$ having supports in the ball $B_N^n$
and with parameters of constancy $\ge l$~\cite[VI.2.]{Vl-V-Z}.
The following embedding holds:
${\cD}^l_N(\bQ_p^n)\subset {\cD}^{l'}_{N'}(\bQ_p^n)$, \ $N\le N'$,
$l\ge l'$. Thus ${\cD}(\bQ_p^n)
=\lim{\rm ind}_{N\to \infty}\lim{\rm ind}_{l\to -\infty}{\cD}^l_N(\bQ_p^n)$.
The space ${\cD}(\bQ_p^n)$ is a complete locally convex vector space.

According to~\cite[VI,(5.2')]{Vl-V-Z}, any function
$\varphi \in {\cD}^l_N(\bQ_p^n)$ is represented in the
following form
\begin{equation}
\label{9.4}
\varphi(x)=\sum_{\nu=1}^{p^{n(N-l)}}
\varphi(c^{\nu})\Delta_{l}(x-c^{\nu}), \quad x\in \bQ_p^n,
\end{equation}
where $\Delta_{l}(x-c^{\nu})$ are the characteristic
functions of the disjoint balls $B_{l}(c^{\nu})$, and the points
$c^{\nu}=(c_1^{\nu},\dots c_n^{\nu})\in B_N^n$ do not depend on $\varphi$.

Denote by ${\cD}'(\bQ_p^n)$ the set of all linear functionals on
${\cD}(\bQ_p^n)$~\cite[VI.3.]{Vl-V-Z}.

Let us introduce in ${\cD}(\bQ_p^n)$ a {\it canonical
$\delta$-sequence} $\delta_k(x)=p^{nk}\Omega(p^k|x|_p)$,
and a {\it canonical $1$-sequence}
$\Delta_k(x)=\Omega(p^{-k}|x|_p)$, $k \in \bZ$, \
$x\in \bQ_p^n$, where
\begin{equation}
\label{10}
\Omega(t)=\left\{
\begin{array}{lcr}
1, &&\quad 0 \le t \le 1, \\
0, &&\quad t>1. \\
\end{array}
\right.
\end{equation}
Here $\Delta_k(x)$ is the characteristic function of the ball $B_{k}^n$.
It is clear~\cite[VI.3., VII.1.]{Vl-V-Z} that
$\delta_k \to \delta$, $k \to \infty$ in ${\cD}'(\bQ_p^n)$
and $\Delta_k \to 1$, $k \to \infty$ in ${\cE}(\bQ_p^n)$.

The Fourier transform of $\varphi\in {\cD}(\bQ_p^n)$ is defined by the
formula
$$
F[\varphi](\xi)=\int_{\bQ_p^n}\chi_p(\xi\cdot x)\varphi(x)\,d^nx,
\quad \xi \in \bQ_p^n,
$$
where $\chi_p(\xi\cdot x)=\chi_p(\xi_1x_1)\cdots\chi_p(\xi_nx_n)
=e^{2\pi i\sum_{j=1}^{n}\{\xi_j x_j\}_p}$; \
$\xi\cdot x$ is the scalar product of vectors.

The Fourier transform is a linear isomorphism ${\cD}(\bQ_p^n)$ into
${\cD}(\bQ_p^n)$. Moreover, according
to~\cite[Lemma~A.]{Taib1},~\cite[III,(3.2)]{Taib3},
~\cite[VII.2.]{Vl-V-Z},
\begin{equation}
\label{12}
\varphi(x) \in {\cD}^l_N(\bQ_p^n) \quad \text{iff} \quad
F\big[\varphi(x)\big](\xi) \in {\cD}^{-N}_{-l}(\bQ_p^n).
\end{equation}

We define the Fourier transform $F[f]$ of a distribution
$f\in {\cD}'(\bQ_p^n)$ by the relation~\cite[VII.3.]{Vl-V-Z}:
\begin{equation}
\label{13}
\langle F[f],\varphi\rangle=\langle f,F[\varphi]\rangle,
\quad \forall \, \varphi\in {\cD}(\bQ_p^n).
\end{equation}

Let $A$ be a matrix and $b\in \bQ_p^n$. Then for a distribution
$f\in{\cD}'(\bQ_p^n)$ the following relation holds~\cite[VII,(3.3)]{Vl-V-Z}:
\begin{equation}
\label{14}
F[f(Ax+b)](\xi)
=|\det{A}|_p^{-1}\chi_p\big(-A^{-1}b\cdot \xi\big)F[f(x)]\big(A^{-1}\xi\big),
\end{equation}
where $\det{A} \ne 0$.
According to~\cite[IV,(3.1)]{Vl-V-Z},
\begin{equation}
\label{14.1}
F[\Delta_{k}](x)=\delta_{k}(x), \quad k\in \bZ, \qquad x \in \bQ_p^n.
\end{equation}
In particular, $F[\Omega(|\xi|_p)](x)=\Omega(|x|_p)$.

The convolution $f*g$ for distributions $f,g\in{\cD}'(\bQ_p^n)$ is
defined (see~\cite[VII.1.]{Vl-V-Z}) as
\begin{equation}
\label{11}
\langle f*g,\varphi\rangle
=\lim_{k\to \infty}\langle f(x)\times g(y),\Delta_k(x)\varphi(x+y)\rangle
\end{equation}
if the limit exists for all $\varphi\in {\cD}(\bQ_p^n)$, where
$f(x)\times g(y)$ is the direct product of distributions. If for
distributions $f,g\in {\cD}'(\bQ_p^n)$ the convolution $f*g$ exists
then~\cite[VII,(5.4)]{Vl-V-Z}
\begin{equation}
\label{15}
F[f*g]=F[f]F[g].
\end{equation}

\begin{Definition}
\label{de0} \rm
Let $\pi_{\alpha}$ be a multiplicative character of the
field $\bQ_p$.
A distribution $f \in {\cD}'(\bQ_p^n)$ is called {\it homogeneous} of
degree $\pi_{\alpha}$ if for all $\varphi \in {\cD}(\bQ_p^n)$ and
$t \in \bQ_p^*$ we have the relation
$$
\Bigl\langle f,\varphi\Big(\frac{x_1}{t},\dots,
\frac{x_n}{t}\Big)\Bigr\rangle
=\pi_{\alpha}(t)|t|_p^n\bigl\langle f,\varphi(x_1,\dots,x_n) \bigr\rangle
$$
i.e., $f(tx)=f(tx_1,\dots,tx_n)=\pi_{\alpha}(t)f(x)$,
$x=(x_1,\dots,x_n)\in \bQ_p^{n}$.
A {\it homogeneous} distribution of degree $\pi_{\alpha}(t)=|t|_p^{\alpha-1}$
($\alpha \ne 0$) is called homogeneous of degree~$\alpha-1$.
\end{Definition}

\section{The $p$-adic Lizorkin spaces}
\label{s3}

Let us introduce the $p$-adic {\it Lizorkin space of test functions\/}
$$
\Phi(\bQ_p^n)=\{\phi: \phi=F[\psi], \, \psi\in \Psi(\bQ_p^n)\},
$$
where
$$
\Psi(\bQ_p^n)=\{\psi(\xi)\in \cD(\bQ_p^n): \psi(0)=0\}.
$$
Here $\Psi(\bQ_p^n), \Phi(\bQ_p^n)\subset \cD(\bQ_p^n)$.
The space $\Phi(\bQ_p^n)$ is called the $p$-adic {\it Lizorkin space of
test functions\/}. The space $\Phi(\bQ_p^n)$ can be equipped with the
topology of the space $\cD(\bQ_p^n)$ which makes $\Phi$ a complete space.

In view of (\ref{12}), the following lemma holds.
\begin{Lemma}
\label{lem1}
{\rm (~\cite{Al-Kh-Sh3},~\cite{Al-Kh-Sh4})}
{\rm (a)} $\phi\in \Phi(\bQ_p^n)$ iff $\phi\in \cD(\bQ_p^n)$ and
\begin{equation}
\label{54}
\int_{\bQ_p^n}\phi(x)\,d^nx=0.
\end{equation}

{\rm (b)} $\phi \in {\cD}^l_N(\bQ_p^n)\cap\Phi(\bQ_p^n)$, i.e.,
$\int_{B^n_{N}}\phi(x)\,d^nx=0$,
iff \ $\psi=F^{-1}[\phi]\in {\cD}^{-N}_{-l}(\bQ_p^n)\cap\Psi(\bQ_p^n)$,
i.e., $\psi(\xi)=0$, $\xi \in B^n_{-N}$.
\end{Lemma}

Unlike the classical Lizorkin space, any function
$\psi(\xi)\in \Phi(\bQ_p^n)$ is equal to zero not only at $\xi=0$
but in a ball $B^n \ni 0$, as well.

Let $\Phi'(\bQ_p^n)$ denote the topological dual of the space $\Phi(\bQ_p^n)$.
We call it the $p$-adic {\it Lizorkin space of distributions\/}.

By $\Psi^{\perp}$ and $\Phi^{\perp}$ we denote the
subspaces of functionals in $\cD'(\bQ_p^n)$ orthogonal to $\Psi(\bQ_p^n)$ and
$\Phi(\bQ_p^n)$, respectively. Thus
$\Psi^{\perp}=\{f\in \cD'(\bQ_p^n): f=C\delta, \, C\in \bC\}$ and
$\Phi^{\perp}=\{f\in \cD'(\bQ_p^n): f=C, \, C\in \bC\}$.

\begin{Proposition}
\label{pr2}
{\rm (~\cite{Al-Kh-Sh3})}
$$
\Phi'(\bQ_p^n)=\cD'(\bQ_p^n)/\Phi^{\perp}, \qquad
\Psi'(\bQ_p^n)=\cD'(\bQ_p^n)/\Psi^{\perp}.
$$
\end{Proposition}

The space $\Phi'(\bQ_p^n)$ can be obtained from $\cD'(\bQ_p^n)$ by
``sifting out'' constants. Thus two distributions in $\cD'(\bQ_p^n)$
differing by a constant are indistinguishable as elements of $\Phi'(\bQ_p^n)$.

Similarly to (\ref{13}), we define the Fourier transform of
distributions $f\in \Phi_{\times}'(\bQ_p^n)$ and $g\in \Psi_{\times}'(\bQ_p^n)$
by the relations:
\begin{equation}
\label{51}
\begin{array}{rcl}
\displaystyle
\langle F[f],\psi\rangle=\langle f,F[\psi]\rangle,
&& \forall \, \psi\in \Psi(\bQ_p^n), \medskip \\
\displaystyle
\langle F[g],\phi\rangle=\langle g,F[\phi]\rangle,
&& \forall \, \phi\in \Phi(\bQ_p^n). \\
\end{array}
\end{equation}
By definition, $F[\Phi(\bQ_p^n)]=\Psi(\bQ_p^n)$ and
$F[\Psi(\bQ_p^n)]=\Phi(\bQ_p^n)$, i.e., (\ref{51})
give well defined objects.

\section{Construction of multiresolution analysis}
\label{s4}

\subsection{$p$-Adic multiresolution analysis.}\label{s4.1}
Denote the factor group $\bQ_p/\bZ_p$ by $I_p$, i.e.
$$
I_p=\{a=p^{-\gamma}\big(a_{0}+a_{1}p+\cdots+a_{\gamma-1}p^{\gamma-1}\big):
\qquad\qquad\qquad\qquad
$$
\begin{equation}
\label{62.0**}
\qquad\qquad
\gamma\in \bN; \, a_j=0,1,\dots,p-1; \, j=0,1,\dots,\gamma-1\}.
\end{equation}

It is well known that $\bQ_p=B_{0}\cup\cup_{\gamma=1}^{\infty}S_{\gamma}$,
where $S_{\gamma}=\{x\in \bQ_p: |x|_p = p^{\gamma}\}$.
In view of (\ref{2}), $x\in S_{\gamma}$, $\gamma\ge 1$ if and only if
$x=x_{-\gamma}p^{-\gamma}+x_{-\gamma+1}p^{-\gamma+1}+\cdots+x_{-1}p^{-1}+\xi$,
where $\xi \in B_{0}$. Since $x_{-\gamma}p^{-\gamma}+x_{-\gamma+1}p^{-\gamma+1}
+\cdots+x_{-1}p^{-1}\in I_p$, we have a ``natural'' decomposition of $\bQ_p$
to a union of mutually  disjoint discs:
$$
\bQ_p=\cup_{a\in I_p}B_{0}(a).
$$
So, $I_p$ is a ``natural'' group of shifts for $\bQ_p$.

\begin{Definition}
\label{de1} \rm
A collection of closed spaces
$V_j\subset\cL^2(\bQ_p)$, $j\in\bZ$ is called a
{\it multiresolution analysis {\rm(}MRA{\rm)} in $\cL^2(\bQ_p)$} if the
following axioms hold

(a) $V_j\subset V_{j+1}$ for all $j\in\bZ$;

(b) $\cup_{j\in\bZ}V_j$ is dense in $\cL^2(\bQ_p)$;

(c) $\cap_{j\in\bZ}V_j=\{0\}$;

(d) $f(\cdot)\in V_j \Longleftrightarrow f(p^{-1}\cdot)\in V_{j+1}$ for all $j\in\bZ$;

(e) there a function $\phi \in V_0$
such that the system $\phi(x-a)$, $a\in I_p$, form an orthonormal
basis for $V_0$.
\end{Definition}

The function $\phi$ from axiom (e) is called {\em scaling} or {\em refinable}.
It follows  immediately from axioms (d) and (e)
that the functions $p^{j/2}\phi(p^{-j}\cdot-a)$, $a\in I_p$,
form an orthonormal basis for $V_j$.

According to the standard scheme (see, e.g.,~\cite[\S 1.3]{NPS}) for construction
of MRA-based wavelets, for each $j$, we define a  space $W_j$ ({\em wavelet space})
as the orthogonal complement of $V_j$ in $V_{j+1}$, i.e.,
\begin{equation}
\label{61}
V_{j+1}=V_j\oplus W_j, \qquad j\in \bZ,
\end{equation}
where $W_j\perp V_j$, $j\in \bZ$. It is not difficult to see that
\begin{equation}
\label{61.0}
f\in W_j \Longleftrightarrow f(p^{-1}\cdot)\in W_{j+1},
\quad\text{for all}\quad j\in \bZ
\end{equation}
and $W_j\perp W_k$, $j\ne k$.
Taking into account axioms (b) and (c), we obtain
\begin{equation}
\label{61.1}
{\oplus_{j\in\bZ}W_j}=\cL^2(\bQ_p)
\quad \text{(orthogonal direct sum)}.
\end{equation}

If now we find  a function $\psi \in W_0$
such that the system $\psi(x-a)$, $a\in I_p$, form an orthonormal
basis for $W_0$, then  the system $p^{j/2}\psi(p^{-j}\cdot-a)$, $a\in I_p$,
is an orthonormal basis for $\cL^2(\bQ_p)$.
Such a function $\psi$ is called a {\em wavelet function} and the basis is a
{\em wavelet basis}.

\subsection{$p$-Adic refinement equation.}\label{s4.2}
Let $\phi$ be a refinable function for a MRA. As was mentioned above,
the system $p^{1/2}\phi(p^{-1}\cdot-a)$, $a\in I_p$,
is a basis for $V_1$. It follows from axoim (a) that
\begin{equation}
\label{62.0-2*}
\phi=\sum_{a\in I_p}\alpha_a\phi(p^{-1}\cdot-a),
\quad \alpha_a\in \bC.
\end{equation}
We see that the function $\phi$ is a solution of a special kind of
functional equation. Such equations are called refinement equations.
Investigation of refinement equations and their solutions is the
most difficult part of wavelet theory in real analysis.

A natural way for construction of a MRA (see, e.g.,~\cite[\S 1.2]{NPS}) is the following.
We start with an appropriate function $\phi$ whose integer shifts form an orthonormal
system, and set
$V_0=\overline{{\rm span}\big\{\phi\big(x-a\big):a\in I_p\big\}}$
and
$V_j=\overline{{\rm span}\big\{\phi\big(p^{-j}x-a\big):a\in I_p\big\}}$,
\ $j\in \bZ$.
It is clear that axioms (d) and (e) of Definition~\ref{de1} are fulfilled.

Of course, not any such a function $\phi$ provides axiom $(a)$.
In the {\em real setting}, the relation $V_0\subset V_{1}$ holds
if and only if the refinable function satisfies a refinement equation.
Situation is different in $p$-adics. Generally speaking, a refinement
equation (\ref{62.0-2*}) does not imply the including property
$V_0\subset V_{1}$. Indeed, we need all the functions $\phi(\cdot-b)$,
$b\in I_p$, to belong to the space $V_1$, i.e., the equalities
$\phi(x-b)=\sum_{a\in I_p}\alpha_{a,b}\phi(p^{-1}x-a)$ should
be fulfilled for all $b\in I_p$. Since $p^{-1}b+a$ is not in $I_p$ in
general, we can not state that refinement equation (\ref{62.0-2*}) implies
$\phi(x-b)=\sum_{a\in I_p}\alpha_{a,b}\phi(p^{-1}x-p^{-1}b-a)\in V_1$ for all
$b\in I_p$.

The {\it refinement equation} reflects {\it some ``self-similarity''}.
The structure of the space $\bQ_p$ has a {\it natural} ``self-similarity''
property which is given by formulas (\ref{79.0}), (\ref{79}). By
(\ref{79}), the characteristic function
$\Delta_0(x)=\Omega\big(|x|_p\big)$ of the unit disc $B_{0}$ is
represented as a sum of $p$  characteristic functions of the
disjoint discs $B_{-1}(r)$, $r=0,1,\dots,p-1$, i.e.,
\begin{equation}
\label{62.0-2}
\Delta_0(x)=\sum_{r=0}^{p-1}\Delta_0\Big(\frac{1}{p}x-\frac{r}{p}\Big),
\quad x\in \bQ_p.
\end{equation}
Thus, in $p$-adics, we have a {\it natural} {\it refinement equation}
(\ref{62.0-2*}):
\begin{equation}
\label{62.0-3}
\phi(x)=\sum_{r=0}^{p-1}\phi\Big(\frac{1}{p}x-\frac{r}{p}\Big),
\quad x\in \bQ_p,
\end{equation}
whose solution is $\phi(x)=\Delta_0(x)=\Omega\big(|x|_p\big)$.
This equation is an analog of the {\em refinement equation} generating
Haar MRA in real analysis.

\subsection{Construction of $2$-adic Haar multiresolution analysis.}\label{s4.3}
Now, using the {\it refinement equation} (\ref{62.0-3}) for $p=2$
\begin{equation}
\label{62.0-4}
\phi(x)=\phi\Big(\frac{1}{2}x\Big)+\phi\Big(\frac{1}{2}x-\frac{1}{2}\Big),
\quad x\in \bQ_2,
\end{equation} and its solution, the {\it refinable function}
$\phi(x)=\Delta_0(x)=\Omega\big(|x|_2\big)$,
we construct $2$-adic multiresolution analysis.

Set
\begin{equation}
\label{70}
V_0=\overline{{\rm span}\big\{\phi\big(x-a\big):a\in I_2\big\}},
\end{equation}
\begin{equation}
\label{70.1}
V_j=\overline{{\rm span}\big\{\phi\big(2^{-j}x-a\big):a\in I_2\big\}},
\quad j\in \bZ.
\end{equation}
It is clear that axioms (d) and (e) of Definition~\ref{de1} are fulfilled
and the system $2^{j/2}\phi(2^{-j}\cdot-a)$, $a\in I_p$ is an orthonormal basis
for $V_j$, $j\in \bZ$.

Note that the characteristic function of the unit disc $\Omega\big(|x|_2\big)$
has a wonderful feature: $\Omega(|\cdot+\xi|_2)=\Omega(|\cdot|_2)$, for all
$\xi\in \bZ_2$ because the $p$-adic norm is non-Archimedean.
In particular, $\Omega(|\cdot\pm 1|_2)=\Omega(|\cdot|_2)$, i.e.,
\begin{equation}
\label{62.0-0}
\phi(x\pm 1)=\phi(x), \quad \forall \, x\in \bQ_2.
\end{equation}
Thus $\phi$ is {\em periodic} with the period $1$.

In view of this fact, taking into account that $2^{-1}b+a$ ($\mod 1$)
is in $I_2$, for all $a,b\in I_2$, it follows from the refinement
equation (\ref{62.0-4}) that $V_0\subset V_1$.
By (\ref{70.1}), this yields axiom $(a)$.

Due to the {\it refinement equation} (\ref{62.0-4}), we obtain that
$V_j\subset V_{j+1}$, i.e., the axiom (a) from Definition~\ref{de1} holds.

\begin{Lemma}
\label{lem2}
The axiom $(b)$ of Definition~{\rm\ref{de1}} holds, i.e.,
$\overline{\cup_{j\in\bZ}V_j}=\cL^2(\bQ_2)$.
\end{Lemma}

\begin{proof}
According to (\ref{9.4}), any function $\varphi\in \cD(\bQ_2)$ belongs
to one of the spaces ${\cD}^l_N(\bQ_2)$, and consequently, is represented
in the form
\begin{equation}
\label{71}
\varphi(x)=\sum_{\nu=1}^{p^{N-l}}
\varphi(c^{\nu})\Delta_{l}(x-c^{\nu}), \quad x\in \bQ_2,
\end{equation}
where $\Delta_{l}(\cdot-c^{\nu})$ are the characteristic functions
of the mutually disjoint discs $B_{l}(c^{\nu})\subset \bQ_2$,
$c^{\nu}\in B_{N}$, \ $\nu=1,2,\dots p^{N-l}$; \ $l=l(\varphi)$,
$N=N(\varphi)$. Since $\Delta_{l}(x-c^{\nu})=\Omega(p^{-l}|x-c^{\nu}|_p)
=\Omega(|p^{l}x-p^{l}c^{\nu}|_p)$ and any number $p^{l}c^{\nu}$
can be represented in the form $p^{l}c^{\nu}=a^{\nu}+b^{\nu}$, where
$a^{\nu}\in I_2$, $b^{\nu}\in \bZ_2$, we have
$\Delta_{l}(x-c^{\nu})=\Delta_{l}(x-a^{\nu})$. Thus any function
$\varphi\in \cD(\bQ_2)$ can be represented in the form
\begin{equation}
\label{72}
\varphi(x)=\sum_{\nu=1}^{p^{N-l}}\alpha_{\nu}
\Delta_{l}(x-a^{\nu}), \quad x\in \bQ_2, \quad a^{\nu}\in I_2,
\quad \alpha_{\nu}\in \bC.
\end{equation}
Consequently, on the basis  of
(\ref{70.1}), $\varphi(x)\in V_{-l}$. Thus any test function
$\varphi$ belongs to one of the space $V_{j}$, where $j=j(\varphi)$.

Since the space $\cD(\bQ_2)$ is dense in $\cL^{2}(\bQ_2)$~\cite[VI.2]{Vl-V-Z},
approximating any function from
$\cL^{2}(\bQ_2)$ by test functions (\ref{72}), we prove our assertion.
\end{proof}

\begin{Lemma}
\label{lem3}
The axiom $(c)$ of Definition~{\rm\ref{de1}} holds, i.e.,
$\cap_{j\in\bZ}V_j=\{0\}$.
\end{Lemma}

\begin{proof}
Suppose that $\cap_{j\in\bZ}V_j\ne\{0\}$. Then there exists a
function $f\in V_j$ for all $j\in\bZ$. Hence, due to (\ref{70.1}),
$f(x)=\sum_{a\in I_2}c_{ja}\phi\big(2^{-j}x-a\big)$ for all $j\in\bZ$.

Let $x=2^{-N}(x_{0}+x_{1}2+x_{2}2^2+\cdots)$.
Since $2^{-j}x=2^{-N-j}(x_{0}+x_{1}2+x_{2}2^2+\cdots)$, for all $j\le -N$,
we have $2^{-j}x\in \bZ_2$, and, consequently, $|2^{-j}x-a|_2>1$ for
all $a\in I_2$, $a\ne 0$.
Thus $\phi\big(2^{-j}x-a\big)=0$ for
all $j\le -N$ and $a\in I_2$, $a\ne 0$.
Since $|2^{-j}x|_2\le 1$,
we have $f(x)=c_{j0}$ for all $j\le -N$.
Similarly, for another $x'=2^{-N'}(x_{0}'+x_{1}'2+x_{2}'2^2+\cdots)$,
we have $f(x')=c_{j'0}$ for all $j\le -N'$.
This yields that $f(x)=f(x')$. Consequently, $f(x)\equiv C$, where $C$ is a
constant. However, if $C\ne 0$, $f\not\in \cL^2(\bQ_2)$. Thus,
$C=0$ and the proof of the theorem is complete.
\end{proof}

According to the above scheme, we introduce the space $W_0$ as the orthogonal
complement of $V_0$ in $V_{1}$.

Set
\begin{equation}
\label{73}
\psi^{(0)}(x)=\phi\Big(\frac{1}{2}x\Big)
-\phi\Big(\frac{1}{2}x-\frac{1}{2}\Big).
\end{equation}

\begin{Lemma}
\label{lem3.1}
The shift system  $\psi^{(0)}(x-a)$, $a\in I_2$, is an orthonormal
basis of the space $W_0$.
\end{Lemma}

\begin{proof}
Let us prove that $W_0\perp V_0$.
It follows from (\ref{62.0-4}), (\ref{73}) that
$$
\big(\psi^{(0)}(x-a),\phi(x-b)\big)
=\int_{\bQ_2}\psi^{(0)}(x-a)\phi(x-b)\,dx
\qquad\qquad\qquad\qquad\qquad
$$
$$
=\int_{\bQ_2}\bigg(\phi\Big(\frac{x}{2}-\frac{a}{2}\Big)
-\phi\Big(\frac{x}{2}-\frac{1}{2}-\frac{a}{2}\Big)\bigg)
\bigg(\phi\Big(\frac{x}{2}-\frac{b}{2}\Big)
+\phi\Big(\frac{x}{2}-\frac{1}{2}-\frac{b}{2}\Big)\bigg)\,dx
$$
for all $a,b\in I_2$. Let $a\ne b$. Since it is impossible $a\ne b+1$,
$b\ne a+1$, taking into account that the functions
$2^{1/2}\phi(2^{-1}\cdot-c)$, $c\in I_2$ are
orthonormal, we obtain $\big(\psi^{(0)}(x-a),\phi(x-b)\big)=0$.
If $a=b$, again due to the orthonormality of the system
$2^{1/2}\phi(2^{-1}\cdot-c)$, $c\in I_2$, taking into account
that $\frac{a}{2},\frac{a}{2}+\frac{1}{2}\in I_2$, we have
$$
\big(\psi^{(0)}(x-a),\phi(x-a)\big)
=\int_{\bQ_2}\bigg(\phi^2\Big(\frac{x}{2}-\frac{a}{2}\Big)
-\phi^2\Big(\frac{x}{2}-\frac{1}{2}-\frac{a}{2}\Big)\bigg)\,dx
\qquad\qquad\quad
$$
$$
\qquad\qquad\qquad\qquad
=\int_{\bQ_2}\phi\Big(\frac{x}{2}-\frac{a}{2}\Big)\,dx
-\int_{\bQ_2}\phi\Big(\frac{x}{2}-\frac{1}{2}-\frac{a}{2}\Big)\,dx=0.
$$
Thus, $\psi^{(0)}(x+a)\perp\phi(x+b)$ for  all $a,b\in I_2$.

The refinement equation
(\ref{62.0-4}) and relation (\ref{73}) imply that
$$
\phi\Big(\frac{x}{2}-a\Big)
=\frac{1}{2}\Big(\phi\big(x-2a\big)+\psi^{(0)}\big(x-2a\big)\Big),
\quad a\in I_2.
$$
Since $\{2^{1/2}\phi(2^{-1}x-a):a\in I_2\}$ is a basis for $V_1$, we have
$V_{1}=V_0\oplus W_0$, i.e., (\ref{61}) holds.
\end{proof}

Thus we prove that the collection $\{V_j:j\in\bZ\}$ is a MRA in
${\cL}^2(\bQ_2)$ and the function $\psi^{(0)}$ defined by (\ref{73})
is a wavelet function. This MRA is a $2$-adic analog of the real
Haar MRA and the wavelet basis generated by $\psi^{(0)}$ is an analog
of real Haar wavelet basis. But in contrast to the real setting,
the {\it refinable function} $\phi$ generating our Haar MRA is
{\em periodic} with the period $1$ (see (\ref{62.0-0})), which
{\em never holds} for real refinable functions. It will be shown
bellow that due of this specific property of $\phi$, there exist
infinity many different orthonormal wavelet bases in the same Haar
MRA (see Sec.~\ref{s5}).

Due to (\ref{8.2-1**}), (\ref{79}), the function $\psi^{(0)}$ can be
rewritten in the form $\psi^{(0)}(x)=\chi_2(2^{-1}x)\Omega(|x|_2)$ and
the Haar wavelet basis is
$$
\psi^{(0)}_{\gamma a}(x)=2^{-\gamma/2}\psi^{(0)}(2^{\gamma}x-a)
\qquad\qquad\qquad\qquad\qquad\qquad\qquad\qquad\qquad
$$
\begin{equation}
\label{62.0-7}
=2^{-\gamma/2}\chi_2\big(2^{-1}(2^{\gamma}x-a)\big)
\Omega\big(|2^{\gamma}x-a|_2\big), \quad x\in \bQ_2,
\quad \gamma\in \bZ, \quad a\in I_2.
\end{equation}

It is clear that
\begin{equation}
\label{62.1-1}
\int_{\bQ_2}\psi^{(0)}_{\gamma a}(x)\,dx=0,
\end{equation}
and, according to Lemma~\ref{lem1}, $\psi^{(0)}_{\gamma a}(x)$
belongs to the Lyzorkin space $\Phi(\bQ_2)$.

\begin{Remark}
\label{rem1} \rm
The Haar wavelet basis (\ref{62.0-7}) coincides with Kozyrev's wavelet basis
(\ref{62.0-1}) for the case $p=2$.
In present paper we restrict ourself by constructing the Haar wavelets
only for  $p=2$. Since Haar refinement equation (\ref{62.0-3}) was
presented for all $p$, a similar  construction may be easily  realized
in the general case.
Moreover, it is not difficult to see that Kozytev's wavelet function
$\theta_{j}(x)$ from (\ref{62.0-1}) can be expressed in terms of the
{\it refinable function} $\phi(x)$ as
\begin{equation}
\label{62.0}
\theta_{j}(x)=\chi_p(p^{-1}jx)\Omega\big(|x|_p\big)
=p^{-1/2}\sum_{r=0}^{p-1}h_r\phi\Big(\frac{1}{p}x-\frac{r}{p}\Big),
\quad x\in \bQ_p,
\end{equation}
where $h_r=p^{1/2}e^{2\pi
i\{\frac{jr}{p}\}_p}$, $r=0,1,\dots,p-1$, \ $j=1,2,\dots,p-1$.
\end{Remark}

\begin{Remark}
\label{rem2} \rm
In view of periodicity (\ref{62.0-0}) of the refinable function
$\phi$, one can use shifts  $\psi^{(0)}(\cdot+a)$, $a\in I_2$,
instead of shifts $\psi^{(0)}(\cdot-a)$, $a\in I_2$.
\end{Remark}

Now we show that there is another function $\psi^{(1)}(x)$ whose shifts
form an orthonormal basis in $W_0$.
Indeed, taking into account (\ref{62.0-0}), we have
$$
\psi^{(1)}(x)=\frac{1}{\sqrt{2}}\bigg(\phi\Big(\frac{x}{2}\Big)
-\phi\Big(\frac{x}{2}-\frac{1}{2}\Big)
-\phi\Big(\frac{x}{2}+\frac{1}{2^2}\Big)
+\phi\Big(\frac{x}{2}-\frac{1}{2^2}\Big)\bigg)
$$
\begin{equation}
\label{73-1}
=\frac{1}{\sqrt{2}}\bigg(\phi\Big(\frac{x}{2}\Big)
-\phi\Big(\frac{x}{2}-\frac{1}{2}\Big)
-\phi\Big(\frac{x}{2}-\frac{1}{2^2}-\frac{1}{2}\Big)
+\phi\Big(\frac{x}{2}-\frac{1}{2^2}\Big)\bigg)
\end{equation}
an its shifts
$$
\psi^{(1)}\Big(x+\frac{1}{2}\Big)
=\frac{1}{\sqrt{2}}\bigg(\phi\Big(\frac{x}{2}+\frac{1}{2^2}\Big)
-\phi\Big(\frac{x}{2}-\frac{1}{2^2}\Big)
-\phi\Big(\frac{x}{2}+\frac{1}{2}\Big)
+\phi\Big(\frac{x}{2}\Big)\bigg)
$$
\begin{equation}
\label{73-1.1}
=\frac{1}{\sqrt{2}}\bigg(\phi\Big(\frac{x}{2}-\frac{1}{2^2}-\frac{1}{2}\Big)
-\phi\Big(\frac{x}{2}-\frac{1}{2^2}\Big)
-\phi\Big(\frac{x}{2}-\frac{1}{2}\Big)
+\phi\Big(\frac{x}{2}\Big)\bigg),
\end{equation}
$$
\psi^{(1)}\big(x-a\big)
\quad\qquad\qquad\qquad\quad\qquad\qquad\qquad\quad\qquad\qquad
\qquad\qquad\quad
$$
$$
=\frac{1}{\sqrt{2}}\bigg(\phi\Big(\frac{x}{2}-\frac{a}{2}\Big)
-\phi\Big(\frac{x}{2}-\frac{a}{2}-\frac{1}{2}\Big)
-\phi\Big(\frac{x}{2}-\frac{a}{2}+\frac{1}{2^2}\Big)
+\phi\Big(\frac{x}{2}-\frac{a}{2}-\frac{1}{2^2}\Big)\bigg).
$$
\begin{equation}
\label{73-1.2}
=\frac{1}{\sqrt{2}}\bigg(\phi\Big(\frac{x}{2}-\frac{a}{2}\Big)
-\phi\Big(\frac{x}{2}-\frac{a}{2}-\frac{1}{2}\Big)
-\phi\Big(\frac{x}{2}-\frac{a}{2}-\frac{1}{2^2}-\frac{1}{2}\Big)
+\phi\Big(\frac{x}{2}-\frac{a}{2}-\frac{1}{2^2}\Big)\bigg).
\end{equation}

Since the system of functions $\{\phi(2^{-1}x-a):a\in I_2\}$ is
orthonormal, in view of (\ref{62.0-0}), formulas (\ref{73-1})--(\ref{73-1.2})
imply that the function $\psi^{(1)}(x)$ and the function $\psi^{(1)}(x-a)$
are orthonormal, whenever $a\in I_2$, $a\ne 0, \frac{1}{2}$. Here we take
into account that all shifts (up to $\mod 1$) of refinable function in
(\ref{73-1}), (\ref{73-1.2}) are distinct.

Similarly, by (\ref{73-1}), (\ref{73-1.1}), we have
$$
\big(\psi^{(1)}(x),\psi^{(1)}(x+2^{-1})\big)
=\int_{\bQ_2}\psi^{(1)}(x)\psi^{(1)}(x+2^{-1})\,dx
\qquad\qquad\qquad
$$
$$
=2^{-1}\int_{\bQ_2}
\bigg\{\phi^2\Big(\frac{x}{2}\Big)
+\phi^2\Big(\frac{x}{2}-\frac{1}{2}\Big)
\qquad\qquad\qquad\qquad\qquad
$$
$$
\qquad
-\phi^2\Big(\frac{x}{2}-\frac{1}{2^2}-\frac{1}{2}\Big)
-\phi^2\Big(\frac{x}{2}-\frac{1}{2^2}\Big)\bigg\}\,dx=0.
$$
and
$$
\big(\psi^{(1)}(x),\psi^{(1)}(x)\big)
=2^{-1}\int_{\bQ_2}\bigg(\phi^2\Big(\frac{x}{2}\Big)
+\phi^2\Big(\frac{x}{2}-\frac{1}{2}\Big)
\qquad\qquad\qquad\qquad
$$
$$
\qquad\qquad
+\phi^2\Big(\frac{x}{2}-\frac{1}{2^2}-\frac{1}{2}\Big)
+\phi^2\Big(\frac{x}{2}-\frac{1}{2^2}\Big)\bigg)\,dx=1.
$$

Thus all shifts of $\psi^{(1)}$ are orthonormal.

It is clear that the functions (\ref{73-1}) and (\ref{73-1.1}) can be
rewritten in the form
\begin{equation}
\label{102-1}
\psi^{(1)}(x)=\frac{1}{\sqrt{2}}
\Big(\psi^{(0)}\big(x\big)-\psi^{(0)}\Big(x+\frac{1}{2}\Big)\Big),
\end{equation}
$$
\psi^{(1)}\Big(x+\frac{1}{2}\Big)=\frac{1}{\sqrt{2}}
\Big(\psi^{(0)}\big(x\big)+\psi^{(0)}\Big(x+\frac{1}{2}\Big)\Big).
$$
It follows that
$$
\psi^{(0)}(x)=\frac{1}{\sqrt{2}}
\Big(\psi^{(1)}\big(x\big)+\psi^{(1)}\Big(x+\frac{1}{2}\Big)\Big).
$$
Since the system $\psi^{(0)}(\cdot-a)$, $a\in I_2$, forms an
orthonormal basis for  $W_0$, the system $\psi^{(1)}(\cdot-a)$,
$a\in I_2$, is another orthonormal basis for $W_0$.

So, we showed that a wavelet basis generated by the Haar MRA is not unique.

\section{Description of $2$-adic Haar bases}
\label{s5}

\subsection{Complex wavelets.}\label{s5.1}
Using the fact that all dilatations and shifts ($x\to 2^{\gamma}x+a$,
$a\in I_2$) of the Haar wavelet function
$\psi^{(0)}$ form a orthonormal basis in ${\cL}^2(\bQ_2)$, we
show that there exist infinitely many wavelet functions $\psi^{(s)}$,
$s\in \bN$ in $W_0$.

In what follows, we shall write the $2$-adic number
$a=2^{-s}\big(a_{0}+a_{1}2+\cdots+a_{s-1}2^{s-1}\big)\in I_2$, \
$a_{j}=0,1$, \ $j=0,1,\dots,s-1$ briefly as a rational number
$a=\frac{m}{2^s}$, where $m=a_{0}+a_{1}2+\cdots+a_{s-1}2^{s-1}$.

Since the characteristic function of the unit disc
$\phi(x)=\Delta_0(x)=\Omega\big(|x|_2\big)$ is {\it periodic} with
the period $\xi \in S_{0}$, the wavelet function $\psi^{0}(x)$ has
the following evident and important property:
\begin{equation}
\label{100}
\psi^{(0)}(x+\xi)=-\psi^{(0)}(x), \qquad \xi \in S_{0}.
\end{equation}
Here $\xi=1+\xi_{1}2+\xi_{2}2^{2}+\cdots$, where $\xi_{j}=0,1$; $j\in \bN$.

Before we prove a general result, we consider the simplest particular case.
Consider the function
\begin{equation}
\label{102}
\psi^{(1)}(x)=\alpha_{0}\psi^{(0)}\big(x\big)
+\alpha_{1}\psi^{(0)}\Big(x+\frac{1}{2}\Big),
\quad \alpha_{0},\alpha_{1}\in \bC,
\end{equation}
and solve the problem when all shifts of this function generates
an orthonormal basis $\psi^{(1)}(x+a)$, $a\in I_2$ in $W_0$.

Taking into account orthonormality of the system $\psi^{(0)}(x+a)$,
$a\in I_2$ and relation (\ref{100}), we can see that the function
$\psi^{(1)}(x)$ and the functions $\psi^{(1)}(x+a)$ are orthonormal
for all $a\in I_2$, $a\ne 0, \frac{1}{2}$. Thus, in view of (\ref{100}),
the system of functions $\psi^{(1)}(x+a)$, $a\in I_2$ is orthonormal if
and only if the system of functions (\ref{102}) and
\begin{equation}
\label{103}
\psi^{(1)}\Big(x+\frac{1}{2}\Big)=-\alpha_{1}\psi^{(0)}\big(x\big)
+\alpha_{0}\psi^{(0)}\Big(x+\frac{1}{2}\Big)
\end{equation}
is orthonormal. Hence, we have $|\alpha_{0}|^2+|\alpha_{1}|^2=1$.
In other words, the matrix
$$
D=\left(
\begin{array}{cc}
\alpha_{0} & \alpha_{1} \\
-\alpha_{1} & \alpha_{0} \\
\end{array}
\right)
$$
is unitary. Thus, the function (\ref{102}), where
$|\alpha_{0}|^2+|\alpha_{1}|^2=1$ is the wavelet function.
It is clear that the wavelet function (\ref{102-1}) is a particular
case of the wavelet function (\ref{102}).

Consequently, all dilatations and shifts of $\psi^{(1)}(x)$
form $2$-adic orthonormal wavelet basis in ${\cL}^2(\bQ_2)$.

Now we will prove a general theorem.

\begin{Theorem}
\label{th4}
Let $s=1,2,\dots$. The function
\begin{equation}
\label{101}
\psi^{(s)}(x)=\sum_{k=0}^{2^s-1}\alpha_{k}\psi^{(0)}\Big(x+\frac{k}{2^s}\Big),
\end{equation}
is the wavelet function (whose dilatations and shifts form $2$-adic
orthonormal wavelet basis in ${\cL}^2(\bQ_2)$) if and only if
\begin{equation}
\label{108}
\alpha_{k}
=2^{-s}(-1)^{k}\sum_{r=0}^{2^s-1}\gamma_re^{-i\pi\frac{2r+1}{2^{s}}k},
\quad k=0,1,2,\dots,2^s-1,
\end{equation}
$\gamma_k\in \bC$, $|\gamma_k|=1$.
\end{Theorem}

\begin{proof}
Suppose that $\psi^{(s)}(x)$, $s\ge 1$ is given by formula (\ref{101}).
Since the system $\psi^{(0)}(\cdot+a)$, $a\in I_2$ is orthonormal
(see Subsec.~\ref{s4.3}) and in view of relation (\ref{100}), it is
easy to see that $\psi^{(s)}$ and $\psi^{(s)}(\cdot+a)$ are orthonormal
for any $a\in I_2$, $a\ne \frac{k}{2^s}$, $k=0,1,\dots 2^s-1$. Thus
the system of functions $\psi^{(s)}(x+a)$, $a\in I_2$ is orthonormal
if and only if the system of functions, consisting of the function
(\ref{101}) and its shifts, i.e.,
$$
\psi^{(s)}\big(x+\frac{r}{2^s}\big)=
-\alpha_{2^s-r}\psi^{(0)}(x)
-\alpha_{2^s-r+1}\psi^{(0)}\Big(x+\frac{1}{2^s}\Big)-\cdots
-\alpha_{2^s-1}\psi^{(0)}\Big(x+\frac{r-1}{2^s}\Big)
$$
\begin{equation}
\label{104}
\qquad
+\alpha_{0}\psi^{(0)}\Big(x+\frac{r}{2^s}\Big)+\cdots
+\alpha_{2^s-r-1}\psi^{(0)}\Big(x+\frac{2^s-1}{2^s}\Big),
\end{equation}
$r=0,1,\dots,2^s-1$ is orthonormal.

Set $\Xi^{(0)}=\{\psi^{(0)}(\cdot+\frac{k}{2^s}):k=0,1,\dots,2^s-1\}^{T}$,
\ $\Xi^{(s)}=\{\psi^{(s)}(\cdot+\frac{k}{2^s}):k=0,1,\dots,2^s-1\}^{T}$.
In view of (\ref{101}), (\ref{104}), $\Xi^{(s)}=D\Xi^{(0)}$,
where
\begin{equation}
D=\label{105}
\left(
\begin{array}{cccccc}
\alpha_{0}&\alpha_{1}&\alpha_{2}&\ldots&\alpha_{2^s-2}&\alpha_{2^s-1} \\
-\alpha_{2^s-1}&\alpha_{0}&\alpha_{1}&\ldots&\alpha_{2^s-3}&\alpha_{2^s-2} \\
-\alpha_{2^s-2}&-\alpha_{2^s-1}&\alpha_{0}&\ldots&\alpha_{2^s-4}&
\alpha_{2^s-3} \\
\hdotsfor{6} \\
-\alpha_{2}&-\alpha_{3}&-\alpha_{4}&\ldots&\alpha_{0}&\alpha_{1} \\
-\alpha_{1}&-\alpha_{2}&-\alpha_{3}&\ldots&-\alpha_{2^s-1}&\alpha_{0} \\
\end{array}
\right).
\end{equation}
Thus the system $\Xi^{(s)}$ is orthonormal if and only if the matrix
$D$ is unitary.

Let $u=(\alpha_{0},\alpha_{1},\dots,\alpha_{2^s-1})^{T}$ be a vector
and
$$
A=\left(
\begin{array}{cccccc}
0&0&\ldots&0&0&-1 \\
1&0&\ldots&0&0&0 \\
0&1&\ldots&0&0&0 \\
\hdotsfor{6} \\
0&0&\ldots&1&0&0 \\
0&0&\ldots&0&1&0 \\
\end{array}
\right).
$$
be a $2^s\times2^s$ matrix.
It is easy to see that
$$
A^ru=(-\alpha_{2^s-r},-\alpha_{2^s-r+1},\dots,-\alpha_{2^s-1},
\alpha_{0},\alpha_{1},\dots,\alpha_{2^s-r-1})^{T},
$$
$r=1,2,\dots,2^s-1$.
Thus $D=\big(u,Au,\dots,A^{2^s-1}u\big)^T$. It is significant that
$A^{2^s}u=-u$. Consequently, in order to describe all matrixes $D$
(or in other words, all vectors $u$), we should find all vectors
$u=(\alpha_{0},\alpha_{1},\dots,\alpha_{2^s-1})^{T}$ such that the
system $\{A^ru:r=0,1,2,\dots,2^s-1\}$ is orthonormal.

In view of the fact that the system $\psi^{(0)}(x+a)$, $a\in I_2$
forms an orthonormal basis in $W_0$, it is easy to see that the
vector $u_0=(1,0,\dots,0,0)^{T}$ is one of mentioned above vectors
$u$. That is the system composed of vectors $u_0$ and
$A^ru_0=(\delta_{0\,r},\delta_{1\,r},\dots,\delta_{2^s-2\,r},
\delta_{2^s-1\,r})^{T}$, $r=1,2,\dots,2^s-1$, is orthonormal,
where $\delta_{i\,r}$ is the Kronecker symbol.

Let us prove that the vector
$u=(\alpha_{0},\alpha_{1},\dots,\alpha_{2^s-1})^{T}$ already
mentioned above such that $A^ru$, \ $r=0,1,2,\dots,2^s-1$ is
orthonormal, can be expressed by the formula $u=Bu_0$ if and only
if $B$ is a unitary matrix such that $AB=BA$.
Indeed, let $u=Bu_0$, where $B$ is a unitary matrix such that
$AB=BA$. Then $A^ru=BA^ru_0$, \ $r=0,1,2,\dots,2^s-1$. Since the
system $A^ru_0$, \ $r=0,1,2,\dots,2^s-1$ is orthonormal and the
matrix $B$ is unitary, the vectors $A^ru$, \ $r=0,1,2,\dots,2^s-1$
are orthonormal.
Conversely, if the system $A^ru$, \ $r=0,1,2,\dots,2^s-1$ is
orthonormal, taking into account that the system $A^ru_0$,
$r=0,1,2,\dots,2^s-1$ is orthonormal, we conclude that
there exists a unitary matrix $B$ such that $A^ru=B(A^ru_0)$, \
$r=0,1,2,\dots,2^s-1$. Since $A^{2^s}u=-u$, $A^{2^s}u_0=-u_0$,
we have an additional relation $A^{2^s}u=BA^{2^s}u_0$. It follows
from the above relations that $(AB-BA)(A^{r}u_0)=0$, \
$r=0,1,2,\dots,2^s-1$. Since the vectors $A^ru_0$, \
$r=0,1,2,\dots,2^s-1$ form a basis in the $2^s$-dimensional
space, we conclude that $AB=BA$.

Thus we have $D=\big(Bu_0,BAu_0,\dots,BA^{2^s-1}u_0\big)^T$.

It is clear that the eigenvalues of $A$ and
the corresponding normalized eigenvectors are
\begin{equation}
\label{106}
\lambda_r=-e^{i\pi\frac{2r+1}{2^{s}}},
\end{equation}
and $v_r=\big((v_r)_1,\dots,(v_r)_{2^s}\big)^T$, respectively,
where
\begin{equation}
\label{107}
(v_r)_l=2^{-s/2}(-1)^{l}e^{-i\pi\frac{2r+1}{2^{s}}l},
\quad l=0,1,2,\dots,2^s-1,
\end{equation}
$r=0,1,2,\dots,2^s-1$. As is well known, the matrix $A$ can be
represented as $A=C\widetilde{A}C^{-1}$, where
$$
\widetilde{A}=\left(
\begin{array}{ccccc}
\lambda_{0}&0&\ldots&0 \\
0&\lambda_{1}&\ldots&0 \\
\vdots&\vdots&\ddots&\vdots \\
0&0&\ldots&\lambda_{2^s-1} \\
\end{array}
\right)
$$
is a diagonal matrix, $C=\big(v_0,v_1,\dots,v_{2^s-1}\big)$. Since $C$
is a unitary matrix, the matrix $B=C\widetilde{B}C^{-1}$ is unitary if
and only if $\widetilde{B}$ is unitary. On the other hand, $AB=BA$ if
and only if $\widetilde{A}\widetilde{B}=\widetilde{B}\widetilde{A}$.
Moreover, since according to (\ref{106}) $\lambda_{k}\ne\lambda_{l}$,
whenever $k\ne l$, all unitary matrix $\widetilde{B}$ such that
$\widetilde{A}\widetilde{B}=\widetilde{B}\widetilde{A}$, are given by
$$
\widetilde{B}=\left(
\begin{array}{ccccc}
\gamma_{0}&0&\ldots&0 \\
0&\gamma_{1}&\ldots&0 \\
\vdots&\vdots&\ddots&\vdots \\
0&0&\ldots&\gamma_{2^s-1} \\
\end{array}
\right),
$$
where $\gamma_k\in \bC$, $|\gamma_k|=1$. Hence, all unitary matrix $B$
such that $AB=BA$, are given by $B=C\widetilde{B}C^{-1}$, where
$\widetilde{B}$ is the above diagonal matrix.

By using formula (\ref{107}), one can calculate
$$
\alpha_{k}=(Bu_0)_k=(C\widetilde{B}C^{-1}u_0)_k
=\sum_{r=0}^{2^s-1}\gamma_r(v_r)_k(\overline{v}_r)_0
\qquad\qquad\qquad\qquad\quad
$$
$$
\qquad\quad
=2^{-s}(-1)^{k}\sum_{r=0}^{2^s-1}\gamma_re^{-i\pi\frac{2r+1}{2^{s}}k},
\quad k=0,1,2,\dots,2^s-1,
$$
where $\gamma_k\in \bC$, $|\gamma_k|=1$. Thus (\ref{108}) holds.

Taking into account that $\Xi^{(0)}=D^{-1}\Xi^{(s)}$, we conclude that
if we define $\psi^{(s)}(x)$ by formula (\ref{101}), where $\alpha_{k}$
is given by (\ref{108}), $k=0,1,2,\dots,2^s-1$, then the system of functions
$\{\psi^{(s)}(\cdot-a):a\in a\in I_2\}$ is orthonormal and forms the
orthonormal basis in $W_0$.

Consequently, all dilatations and shifts of the function (\ref{101})
form $2$-adic orthonormal wavelet basis in ${\cL}^2(\bQ_2)$.
\end{proof}

It is clear that $\int_{\bQ_2}\psi_{\gamma a}^{(s)}(x)\,dx=0$,
and in view of Lemma~\ref{lem1}, $\psi_{\gamma a}^{(s)}(x)$ belongs
to the Lizorkin space $\in \Phi(\bQ_2^n)$.

\subsection{Real wavelets.}\label{s5.2}
Using formulas (\ref{108}), one can extract all {\it real} wavelet
functions (\ref{101}).

Let $s=1$. According to (\ref{102}),  (\ref{103}),
\begin{equation}
\label{109}
\psi^{(1)}(x)=\cos\theta\,\psi^{(0)}\big(x\big)
+\sin\theta\,\psi^{(0)}\Big(x+\frac{1}{2}\Big)
\end{equation}
is the {\it real} wavelet function.

Let $s=2$. Set $\gamma_r=e^{i\theta_r}$, $r=0,1,2,\dots,2^s-1$.
Then (\ref{108}) imply that the wavelet function $\psi^{(1)}(x)$
is real if and only if
$$
\begin{array}{rcl}
\displaystyle
\sin\theta_1+\sin\theta_2+\sin\theta_3+\sin\theta_4&=&0, \\
\displaystyle
\cos\theta_1-\cos\theta_2+\cos\theta_3-\cos\theta_4&=&0, \\
\displaystyle
\sin\theta_1-\sin\theta_2-\sin\theta_3+\sin\theta_4&=& \\
\displaystyle
\cos\theta_1+\cos\theta_2-\cos\theta_3-\cos\theta_4,&& \\
\displaystyle
\sin\theta_1-\sin\theta_2-\sin\theta_3+\sin\theta_4&=& \\
\displaystyle
-(\cos\theta_1+\cos\theta_2-\cos\theta_3-\cos\theta_4).&& \\
\end{array}
$$
The last relations are equivalent to the system
$$
\begin{array}{rclrcl}
\displaystyle
\sin\theta_1&=&-\sin\theta_4, \quad \cos\theta_1&=&\cos\theta_4,\\
\displaystyle
\sin\theta_2&=&-\sin\theta_3, \quad \cos\theta_2&=&\cos\theta_3.\\
\end{array}
$$

Thus for $s=2$ the real wavelet functions (\ref{101}) is represented as
$$
\psi^{(1)}(x)=\frac{1}{2}(\cos\theta_1+\cos\theta_2)\psi^{(0)}\big(x\big)
\qquad\qquad\qquad\qquad\qquad\qquad\qquad
$$
$$
+\frac{1}{2\sqrt{2}}(\cos\theta_1-\cos\theta_2+\sin\theta_1+\sin\theta_2)
\psi^{(0)}\Big(x+\frac{1}{2^2}\Big)
$$
$$
+\frac{1}{2}(\sin\theta_1-\sin\theta_2)\psi^{(0)}\Big(x+\frac{1}{2}\Big)
\qquad\qquad\qquad\qquad
$$
\begin{equation}
\label{110}
+\frac{1}{2\sqrt{2}}(\cos\theta_1-\cos\theta_2-\sin\theta_1-\sin\theta_2)
\psi^{(0)}\Big(x+\frac{1}{2^2}+\frac{1}{2}\Big).
\end{equation}

In particular, for the special cases $\theta_1=\theta_2=\theta$,
$\theta_1=-\theta_2=\theta$, $\theta_1=\theta_2+\frac{\pi}{2}=\theta$,
we obtain one-parameter families of the real wavelet functions
\begin{equation}
\label{111}
\begin{array}{rclrcl}
\displaystyle
\psi^{(1)}(x)&=&\cos\theta\psi^{(0)}\big(x\big)
+\sin\theta\psi^{(0)}\Big(x+\frac{1}{2}\Big),\\
\displaystyle
\psi^{(1)}(x)&=&\cos\theta\psi^{(0)}\big(x\big)
+\frac{1}{\sqrt{2}}\sin\theta\psi^{(0)}\Big(x+\frac{1}{2^2}\Big) \\
\displaystyle
&& \qquad\qquad\qquad
-\frac{1}{\sqrt{2}}\sin\theta\psi^{(0)}\Big(x+\frac{1}{2^2}+\frac{1}{2}\Big),\\
\displaystyle
\psi^{(1)}(x)&=&\frac{1}{2}(\cos\theta-\sin\theta)\psi^{(0)}\big(x\big)
+\frac{1}{2\sqrt{2}}(\cos\theta+\sin\theta)\psi^{(0)}\Big(x+\frac{1}{2^2}\Big)\\
\displaystyle
&& \qquad\qquad\qquad\qquad\qquad
-\frac{1}{2}(\cos\theta-\sin\theta)\psi^{(0)}\Big(x+\frac{1}{2}\Big),\\
\end{array}
\end{equation}
respectively.

\section*{Acknowledgments}

The authors are greatly indebted to E.~Yu.~Panov for fruitful discussions.

\end{document}